\documentclass[10pt,journal]{IEEEtran}
\usepackage{cite}
 \usepackage[cmex10]{amsmath}
\usepackage{amssymb,amsfonts}
\usepackage{algorithmic}
\usepackage{graphicx}
\usepackage{textcomp}
\def\BibTeX{{\rm B\kern-.05em{\sc i\kern-.025em b}\kern-.08em
    T\kern-.1667em\lower.7ex\hbox{E}\kern-.125emX}}

\newtheorem{thm}{\textbf{Theorem}}
\newtheorem{assump}{\textbf{Assumption}}

\begin{document}
\title{Rapid Stabilization of Timoshenko Beam by PDE Backstepping}

\author{Guangwei Chen, Rafael Vazquez, \IEEEmembership{Senior Member, IEEE}, and Miroslav Krstic \IEEEmembership{Fellow, IEEE} 
\thanks{This work was partially supported by grant PGC2018-100680-B-C21 funded by MCIN/AEI/10.13039/501100011033. This journal paper is an expanded version of the conference submission~\cite{chen2022}.}
\thanks{Guangwei Chen is with the Institute of Cyber-Systems and Control, Zhejiang University, Hangzhou, 310027, China (e-mail: gwchen@zju.edu.cn). }
\thanks{Rafael Vazquez is with the Department of Aerospace Engineering, Universidad de Sevilla, Camino de los Descubrimiento s.n., Sevilla, Spain (e-mail: rvazquez1@us.es).}
\thanks{M. Krstic is with
the Department of Mechanical Aerospace Engineering,
 University of California, San Diego, CA 92093-0411, USA 
(e-mail: krstic@ucsd.edu)}}


\maketitle

\begin{abstract}
In this paper, we present a rapid boundary stabilization of a Timoshenko beam with anti-damping and anti-stiffness at the uncontrolled boundary, by using  PDE backstepping. We introduce a transformation to map the Timoshenko beam states into a $(2+2)\times(2+2)$ hyperbolic PIDE-ODE system. Then backstepping is applied to obtain a control law guaranteeing closed-loop stability of the origin in the $H^1$ sense. Arbitrarily rapid stabilization can be achieved by adjusting control parameters. Finally, a numerical simulation shows that the proposed controller can rapidly stabilize the Timoshenko beam. This result extends a previous work which considered a slender Timoshenko beam with Kelvin-Voigt damping, allowing destabilizing boundary conditions at the uncontrolled boundary and attaining an arbitrarily rapid convergence rate.
\end{abstract}


\section{Introduction}
Flexible beams are widely used in many applications ranging from aerospace to civil structures.  Correspondingly, beam stabilization has become an important research topic. Among all the beam models, Timoshenko model, as the most realistic of the 1D distributed parameter models, takes into account both the rotatory inertia of the beam cross-sections and the deflection due to shear effect. In this paper, we focus on the control of such a model by applying the backstepping method.

Extensive literature exists on the control of beams and particularly on Timoshenko beams. We next give an overview of past results. More than three decades ago, Kim and Renardy \cite{Kim1987} used a classical boundary damper feedback which required both 
the space and time derivatives at the tip of the beam. Later, considering a clamped-free Timoshenko beam, Morgul \cite{Mor1992} proposed a more general dynamic 
boundary feedback. Balakrishnan~\cite{bal1997,bal2005} considered boundary conditions leading to superstability (vanishing of the beam states in finite time), for clamped boundary conditions on the uncontrolled end. Taylor and Yau \cite{Taylor2002} studied a rotating
Timoshenko beam that  can be stabilized by both applying a force at the free end
and a torque at the pivoted end. Xu et al. investigated the use of pointwise
feedback controls based on asymptotic analysis of eigenvalues and the eigenfunctions \cite{Xu2003}. Soufyane et al. achieved  uniform stabilization  by using a locally distributed damping; in this case, stability can be guaranteed if and only if the two wave equations have the same speeds \cite{Soufyane2003}. Macchelli et al. used a distributed port Hamiltonian (dpH) approach to describe Timoshenko beams and proposed a finite dimensional passive controller that shapes the beam's total energy\cite{Macchelli2004}. This approach has also been followed by other authors; for instance, Siuka et al. also adopted a dpH model and proposed an invariant-based method to achieve stabilization \cite{Siuka2011} and Wu et al. used a passive LQG control design method~\cite{Wu2017}.  Xu presented a boundary feedback design for the exponential stabilization of  a Timoshenko beam with both ends free, and gave an explicit asymptotic formula
of eigenvalues of the closed loop system \cite{Xu2005}. Considering a Timoshenko beam
with local Kelvin–Voigt damping, Zhao  et al.  obtained exponential stability under some additional hypotheses \cite{Zhao2005}.  Krstic et al. extended the backstepping method, by using a singular perturbation approach, to  controller and observer design for a slender Timoshenko beam, with actuation only at the beam base and sensing only at the beam tip \cite{Krstic2006, Siranosian2009}. For a nonuniform Timoshenko beam with spatial-varying parameters, Ammar-Khodja et al. \cite{Ammar2007} studied the stabilization for both internal and boundary cases with one control force. He et al. designed an output-feedback control law using a Lyapunov-based approach with a disturbance observers \cite{HeZhang2012}; the Lyapunov approach is a powerful tool in design of control laws for beams, not only in the Timoshenko model (see e.g.~\cite{Coron1998}). Extending the approach, He et al. proposed an adaptive integral-Barrier Lyapunov function boundary control for  inhomogeneous Timoshenko beams  with constraints \cite{He2014}. Considering both system uncertainties and uncertain input backlash non-linearity, He et al. gave  vibration boundary control law using a disturbance observer \cite{Wei2015}. Allowing for hysteresis of the boundary control input, Liu and Xu proposed a dynamic feedback control
law that exponentially stabilized the beam with
distributed delay \cite{Liu2013}. Yildirim et al. proposed a novel optimal piezoelectric control approach for suppressing vibrations \cite{Yildirim2016}. Finally, to cite several very recent contributions, Ma et al. introduced a
prescribed performance function restricting within an arbitrarily small residual set \cite{Ma2020}, Guo and Meng consider a two-dimensional robust output tracking for a Timoshenko beam equation by using an observer-based error feedback control approach~\cite{Guo2021}, and Mattioni et al. address a beam clamped on a moving inertia actuated by an external torque and force
with the dpH method using strong dissipation feedback~\cite{mattioni2020}, and also in the case of having a mass at the controlled end~\cite{mattioni2021}.

In recent years, the backstepping method has proven itself as a powerful design method for control of infinite dimensional systems. However,  beyond the results in \cite{Krstic2006,Siranosian2009} more than a decade ago, backstepping has not been fully exploited for Timoshenko beam control, even though it has produced results for the shear beam model~\cite{Krstic2008} and the Euler-Bernouilli model~\cite{Smysh2009,Wang2015}. In the present paper, we aim to achieve  rapid stabilization of a Timoshenko beam with anti-damping and anti-stiffness at the uncontrolled boundary. The decay rate can be prescribed arbitrarily by setting the controller parameters. Specifically, we propose an initial transformations of the Timoshenko beam states to a new set of variables verifying a system of $(2+2)\times(2+2)$  hyperbolic PIDEs, coupled with two ODEs. Then the backstepping method is directly applied to controller design of the new system,  by extending previously-developed tools \cite{c1}. 

Thus, the main contribution of this paper with respect to previous results is allowing destabilizing boundary conditions at the uncontrolled boundary (numerous works consider simple clamped conditions) and attaining an arbitrarily rapid convergence rate.

The paper is organized as follows: Section \ref{sec2} presents
the Timoshenko beam model. Section \ref{sec-control} gives the design of the boundary controller and the main result. Then, Section \ref{contdesign} analyzes the resulting controller. Section \ref{stability} studies the closed-loop stability. Finally, Section \ref{sec-numerical} validates the effectiveness of the proposed controllers by numerical simulation and Section \ref{sec-conclusion} closes the paper with some concluding remarks.

\section{Problem Statement}\label{sec2}
The goal of this work is to design an exponentially stabilizing feedback control law (with arbitrary convergence rate) for the equilibrium at the origin of the following Timoshenko beam model, which is given by the following system of PDEs
\begin{eqnarray}\label{plant_eq1}
\varepsilon {u_{tt}} &=& {u_{xx}} - {\alpha _x},\\
 \mu {\alpha_{tt}}&=& {\alpha _{xx}} + \frac{a}{\varepsilon }\left( {{u_x} - \alpha } \right), \label{plant_eq2}
\end{eqnarray}
where $u(x, t)$ denotes the displacement, $\alpha(x, t)$ denotes the angle of rotation,
for $x\in (0,1)$, $t>0$. The coefficients $\varepsilon,\mu>0$, and $a\in \mathbb{R}$ 
are non-dimensional physical parameters defined in \cite{Han1999}. The boundary conditions of (\ref{plant_eq1})--(\ref{plant_eq2}) are
\begin{eqnarray}\label{plant_bc1}
{u_x}\left( 0, t \right)&=& \alpha (0, t) - \theta {u_t}\left( 0, t \right) - \xi u(0, t),\\
u_x(1, t)&=&V_1(t),\\
\alpha_x(0, t)&=&0,\\
\alpha_x(1, t)&=&V_2(t), \label{plant_bc2}
\end{eqnarray}
with $\theta,\xi \in \mathbb{R}$ (which represent, respectively, anti-damping and anti-stiffness), and $V_1(t)$ and $V_2(t)$ the actuation variables that have to be designed. The initial conditions for the system  (\ref{plant_eq1})--(\ref{plant_bc2}) are denoted by $u_0(x)=u(x,0)$, $\alpha_0(x)=\alpha(x,0)$, $u_{0t}=u_t(x,0)$, $\alpha_{0t}=\alpha_t(x,0)$.
\begin{assump} \label{assump1}
The anti-damping coefficient $\theta$ appearing in (\ref{plant_bc1}) verifies $\theta \neq \ \sqrt \varepsilon$.
\end{assump}

This assumption is critical in the transformation posed next. To understand its underlying reason, consider just a simple wave equation $\varepsilon {u_{tt}} = {u_{xx}} $; then, a boundary condition of the type ${u_x}\left( 0, t \right)= - \sqrt{\epsilon} {u_t}\left( 0, t \right)$ is ill-posed, since a solution by the method of characteristics (or alternatively, d'Alambert's solution) will end up with one over-determined characteristic and one under-determined characteristic; then, depending on the other boundary condition, this leads to only trivial or constant functions solving the equation, which in general will not agree with the initial conditions.

\section{Controller design and main result}\label{sec-control}
This section presents the design of our boundary control law, starting in Section~\ref{sec-transformation} with a transformation of the Timoshenko beam states to a new set of variables verifying a hyperbolic-ODE coupled system. Next, in Section~\ref{main} the control law is introduced and the main result is stated.

\subsection{Transformation of the wave PDE representation of the beam to a system of hyperbolic PDEs coupled with ODEs}\label{sec-transformation}
As a first step, the Timoshenko beam is transformed into a first-order hyperbolic integro-differential system coupled with ODEs. This represents an alternative, novel idea to design a controller for this plant, since it opens the door to apply 1-D hyperbolic system control designs. The system becomes a $(2+2)\times(2+2)$ homodirectional system of hyperbolic PIDEs, coupled with two ODEs, by using  the following transformations
\begin{eqnarray}\label{Tran_eq1}
p&=& {u_x} + \sqrt \varepsilon  {u_t},\\
\label{Tran_eq2}
q &=& {u_x} - \sqrt \varepsilon  {u_t}, \\
\label{Tran_eq3}
r &=&{\alpha _x} + \sqrt \mu  {\alpha _t},\\
\label{Tran_eq4}
s & =& {\alpha _x} - \sqrt \mu  {\alpha _t},\\
x_1 &=& u(0,t),\\
x_2 &=& \alpha(0,t).
\end{eqnarray}
Then (\ref{plant_eq1})--(\ref{plant_bc2}) is equivalent to the PDE-ODE system
\begin{eqnarray}  \label{Equ_eq1}
 {p_t} &=&\frac{1}{\sqrt \varepsilon }{p_x} - \frac{1}{2\sqrt \varepsilon }\left( {r + s } \right),
\\ 
{r_t} &=&\frac{1}{\sqrt \mu }{r_x}+\frac{a}{{2\varepsilon\sqrt \mu   }}\left( {p + q}\right)\nonumber\\
&&-\frac{a}{{2\varepsilon\sqrt \mu   }}\left[\int_0^x {\left( {r\left( y, t \right) + s \left( y, t \right)} \right)dy + 2x_2 } \right],\,\,\,  \label{Equ_eq2}
\\  \label{Equ_eq3}
  {q_t} &=& -\frac{1}{ \sqrt \varepsilon}{q_x} - \frac{1}{2 \sqrt \varepsilon}\left( {r + s } \right),
\\  
{ s_t} &=&-\frac{1}{ \sqrt \mu  }{ s_x}+\frac{a}{{2\varepsilon  \sqrt \mu  }}\left({p + q}\right) \nonumber\\
&&- \frac{a}{{2\varepsilon  \sqrt \mu  }}\left[\int_0^x {\left( {r\left( y, t \right) + s \left( y, t \right)} \right)dy +2x_2} \right],\,\,\, \label{Equ_eq4}
\\  \label{ODE_eq1}
 \dot x_1 &=&\frac{2}{{\sqrt \varepsilon   - \theta }}\left[ {\xi x_1 - x_2 + p(0, t)} \right],
\\  \label{ODE_eq2}
\dot x_2 &=&-\frac{1}{{\sqrt \mu  }}s \left( 0, t \right),
\end{eqnarray}
with boundary conditions
\begin{eqnarray}
\label{BD_eq1}
 q(0, t) &=&-\frac{(\sqrt{\varepsilon}+\theta)}{\sqrt{\varepsilon}-\theta}p(0, t)-\frac{2\sqrt \varepsilon}{\sqrt{\varepsilon}-\theta}(\xi x_1-x_2),\,\,\,
\\  \label{BD_eq2}
 s(0, t)&=&-r(0, t),
\\  \label{BD_eq3}
 p\left( 1, t \right) &=& {V_p}(t),
\\  \label{BD_eq4}
 r\left( 1, t \right) &=&{V_r}(t),
\end{eqnarray}
where  ${V_p}(t)=V_1(t)+\sqrt \varepsilon  {u_t}(1, t)$ and ${V_r}(t)=V_2(t)+ \sqrt \mu  {\alpha _t}(1, t)$ are the redefined control variables for this plant.

It must be noted that (\ref{Equ_eq1})--(\ref{Equ_eq4}) is a system of $(2+2)\times (2+2)$ 1-D hyperbolic PIDEs coupled with two ODEs (\ref{ODE_eq1})--(\ref{ODE_eq2}) that has not been explored before, but facilitates analysis and design of controllers, as hyperbolic systems have been widely explored~\cite{c5}. For instance, it is easy to see that the ``superstability'' (convergence in finite time) result stated in~\cite{bal1997,bal2005}, using static output feedback, is only possible with clamped boundary conditions $u(0,t)=\alpha(0,t)=0$, as they automatically impose that the finite-dimensional states $x_1$ and $x_2$ are zero. Note that in that case a straightforward application of the standard  backstepping method for hyperbolic systems~\cite{c6} can achieve finite-time stabilization in the minimum possible time without the need of the additional conditions in~\cite{bal1997}.

The system (\ref{Equ_eq1})--(\ref{BD_eq4}) is similar to the one stabilized with backstepping in \cite{c1}. Thus, the method therein can be easily adapted. However, as a first step, the procedure requires ordering the states $p$ and $r$ as function of their transport speeds, namely $\frac{1}{\sqrt\varepsilon}$ and $\frac{1}{\sqrt\mu}$. There are three possible cases: $\frac{1}{\sqrt\varepsilon}>\frac{1}{\sqrt\mu}$, $\frac{1}{\sqrt\varepsilon}<\frac{1}{\sqrt\mu}$, and $\frac{1}{\sqrt\varepsilon}=\frac{1}{\sqrt\mu}$. In what follows, we assume $\frac{1}{\sqrt\varepsilon}>\frac{1}{\sqrt\mu}$; the  case $\frac{1}{\sqrt\varepsilon}<\frac{1}{\sqrt\mu}$ can be treated analogously by switching the order of the states $p$ and $r$ in all subsequent steps, and the equality case becomes a straightforward extension of the $2\times2$ backstepping design. Define 
\begin{equation}
Z=\left[                 
  \begin{array}{c}   
    p \\  
    r \\  
  \end{array}
\right], Y=\left[                 
  \begin{array}{c}   
    q \\  
    s \\  
  \end{array}
\right], X=\left[                 
  \begin{array}{c}   
     x_1\\  
   x_2\\ 
  \end{array}
\right], V=\left[                 
  \begin{array}{c}   
     V_p\\  
     V_r\\ 
  \end{array}
\right].
\end{equation}
Then, (\ref{Equ_eq1})--(\ref{BD_eq4}) can be written in the following simplified matrix form
\begin{eqnarray}\label{m_e1}  
Z_t&=&\Sigma Z_x +\Lambda_1 (Z+Y)+\Lambda_2 X\nonumber\\
&&+\int^x_0 F[Z(y,t)+Y(y,t)]dy\\ \label{m_e2}
Y_t&=&-\Sigma Y_x+\Lambda_1 (Y+ Z)+\Lambda_2 X\nonumber\\
&&+\int^x_0 F[Z(y,t)+Y(y,t)]dy
\\ \label{m_e3}
\dot{X}&=&(A+B_2D)X+(B_1+B_2C)Z(0, t)
\end{eqnarray}
with boundary conditions
\begin{eqnarray}
Z(1,t)&=&V\\
Y(0,t)&=&CZ(0, t)+DX\label{BDm_e} 
\end{eqnarray}
where 
\begin{eqnarray} \label{eq4}      
  \Sigma&=&\left[                 
  \begin{array}{cc}   
  \frac{1}{\sqrt\varepsilon}&0 \\  
   0&\frac{1}{\sqrt\mu}\\ 
  \end{array}
\right], ~ \Lambda_1=\left[                 
  \begin{array}{cc}   
    0&-\frac{1}{2\sqrt\epsilon} \\  
    \frac{a}{2\varepsilon\sqrt\mu}&0\\  
  \end{array}
\right],
\\ \label{eq6}  
 ~\Lambda_2&=&\left[                 
  \begin{array}{cc}   
    0&0 \\  
    0&-\frac{a}{\varepsilon\sqrt\mu}\\  
  \end{array}
\right], ~F=\left[                 
  \begin{array}{cc}   
    0&0 \\  
    0&-\frac{a}{2\varepsilon\sqrt\mu}\\  
  \end{array}
\right], \,\,\\       
 A&=&\left[                 
  \begin{array}{cc}   
     \frac{2\xi}{{\sqrt \varepsilon   - \theta }}& -\frac{2}{{\sqrt \varepsilon   - \theta }}\\  
   0&0\\ 
  \end{array}
\right] \quad            
 B_1=\left[                 
  \begin{array}{cc}   
   \frac{2}{{\sqrt \varepsilon   - \theta }}&0 \\  
  0&0\\ 
  \end{array}
\right],\,\,\,\,\,\,\,\,
\\ \label{eq7}    
 B_2&=&\left[                 
  \begin{array}{cc}   
  0&0 \\  
   0&-\frac{1}{\sqrt \mu}\\ 
  \end{array}
\right],~    
 C=\left[                 
  \begin{array}{cc}   
     -\frac{\sqrt{\varepsilon}+\theta}{\sqrt{\varepsilon}-\theta}& 0\\  
   0&-1\\ 
  \end{array}
\right],\\            
 D&=&\left[                 
  \begin{array}{cc}   
  -\frac{2\sqrt{\varepsilon}\xi}{\sqrt{\varepsilon}-\theta}&\frac{2\sqrt{\varepsilon}}{\sqrt{\varepsilon}-\theta} \\  
   0&0\\ 
  \end{array}
\right].
\end{eqnarray}
The system (\ref{m_e1})--(\ref{BDm_e}), differently from \cite{c1}, contains integral coupling terms, and the states of ODEs appearing inside the domain of the PDEs. 
\subsection{Stabilizing control law and main result}\label{main}
For system (\ref{m_e1})--(\ref{BDm_e}), the following control law is obtained in Section~\ref{contdesign}.
\begin{equation}\label{eqn-controlaw}
V= \int_0^1 {{K}\left( {1,y} \right)Z\left( {y,t} \right)}dy+\int_0^1 {{L}\left( {1,y} \right)Y\left( {y,t} \right)}dy+\Phi(1)X,
\end{equation}
whose gain kernels are the particular values of the $2\times2$ matrices 
\begin{eqnarray}
K(x, y)&=&\left[ \begin{array}{cc} 
k_{11}&k_{12}\\
k_{21}&k_{22}\\
\end{array}\right] ,\, L(x, y)=\left[ \begin{array}{cc} 
l_{11}&l_{12}\\
l_{21}&l_{22}\\
\end{array}\right], \,\,\,\,\\ \Phi(x)&=&\left[ \begin{array}{cc} 
\phi_{11}&\phi_{12}\\
\phi_{21}&\phi_{22}\\
\end{array}\right],
\end{eqnarray}
evaluated at $x=1$. These matrices verify the following (well-posed) kernel equations 
\begin{eqnarray}
\Sigma K_{x}+K_{y}\Sigma
&=&
\left(K+L\right)\Lambda_1-\Omega(x)K-F\nonumber\\
&&+\int^x_y \left[K(x,s)+L(x, s)\right]Fds,\label{ker_1}\\
\Sigma L_{x}-L_{y}\Sigma &=&\left(K+L\right)\Lambda_1-\Omega(x)L-F\nonumber\\
&&+\int^x_y \left[K(x,s)+L(x, s)\right]Fds,\label{ker_2}
\end{eqnarray}
\begin{eqnarray}
\Phi_x&=&\Sigma^{-1}\Phi A-\Sigma^{-1}\Lambda_2+\Sigma^{-1}\Phi B_2D\nonumber \\
&&-\Sigma^{-1}\Omega(x)\Phi+\int^x_0 \Sigma^{-1}(K-L)\Lambda_2 dy\nonumber\\
&&+\Sigma^{-1}L(x, 0)\Sigma D,  \label{eqn-Phi}
\end{eqnarray}
with boundary conditions for $K$ and $L$
\begin{eqnarray}
\Sigma L(x,x)+L(x,x)\Sigma&=&-\Lambda_1,\label{ker_4} \\
\Sigma K(x,x)-K(x,x)\Sigma&=&-\Lambda_1+\Omega(x),\label{ker_3}\\
 K(x, 0)-(x, 0)\Sigma C\Sigma^{-1}&=&\Phi B_1\Sigma^{-1}+\Phi B_2C\Sigma^{-1},\,\,\,\,\,\,\,\,\label{ker_5} 
\end{eqnarray}
with
\begin{equation}
\Omega(x)=\left[
\begin{array}{cc}
0&0\\
\omega_{21}&0\\
\end{array}
\right],
\end{equation}
where
$\omega_{21}(x, t)=(\frac{1}{\sqrt\mu}-\frac{1}{\sqrt\varepsilon})k_{21}(x, x)+\frac{a}{2\varepsilon}$
and boundary conditions for $\Phi(X)$ in (\ref{eqn-Phi}) as follows
\begin{equation}\label{Phi_delta}
\Phi(0)=\left[ \begin{array}{cc} 
-\xi-\frac{\delta_1}{\kappa} &1+\frac{1}{\sqrt{\mu}}\\
0&-\delta_2\sqrt{\mu}\\
\end{array}\right],
\end{equation}
where $\kappa=2/(\sqrt{\varepsilon}-\theta)$. The parameters $\delta_1,\delta_2$ are arbitrary values which directly determine the decay rate of the closed-loop controlled Timoshenko beam (see Section~\ref{stability}). The well-posedness  of the kernel equations for $K(x, y), L(x, y), \Phi(x)$ is  stated in Theorem~\ref{thm2} in Section~\ref{transf}.

Expressing (\ref{eqn-controlaw}) in terms of the  Timoshenko beam variables:
\begin{eqnarray}
V_1&=&-\int^1_0 \left( k_{11, y}(1, y)+l_{11, y}(1, y) \right) u(y, t)dy \nonumber \\
&&+\int^1_0\sqrt\varepsilon \left(  k_{11}(1, y)+ l_{11}(1, y \right) u_{t}(y, t)dy\nonumber\\
&&-\int^1_0 
\left( k_{12, y}(1, y)+l_{12, y}(1, y)
\right)
\alpha(y, t)dy \nonumber \\
&&+\int^1_0 
\sqrt\mu  \left( k_{12}(1, y)+l_{12}(1, y \right) \alpha_{t}(y, t)dy\nonumber \\
&&+\left(k_{11}(1,1)+l_{11}(1,1)\right)u(1, t)\nonumber \\
&&
-\left( k_{11}(1, 0)+l_{11}(1, 0)-\phi_{11}(1)\right)u(0, t)
\nonumber \\
&&
-(k_{12}(1, 0)+l_{12}(1, 0)-\phi_{12}(1))\alpha(0, t)
\nonumber \\
&&+(k_{12}(1,1)+l_{12}(1,1))\alpha(1, t)-\sqrt \varepsilon  {u_t}(1, t),\,\,\,\, \label{U_1}
\end{eqnarray}
\begin{eqnarray}
V_2&=&-\int^1_0 (k_{21, y}(1, y)+l_{21, y}(1, y))u(y, t)dy
\nonumber\\
&&
+\int^1_0 \sqrt\varepsilon (k_{21}(1, y)+l_{21}(1, y))u_{t}(y, t)dy\nonumber\\
&&-\int^1_0 (k_{22, y}(1, y)+l_{22, y}(1, y))\alpha(y, t)dy
\nonumber\\
&&
+\int^1_0 \sqrt\mu ( k_{22}(1, y)+l_{22}(1, y)) \alpha_{t}(y, t)dy\nonumber \\
&&+\left( k_{21}(1, 1)+l_{21}(1,1)
\right)u(1, t)
\nonumber\\
&&
-\left(k_{21}(1, 0)+l_{21}(1, 0)-\phi_{21}(1)
\right) u(0, t)\nonumber\\
&&
-\left( k_{22}(1, 0)+l_{22}(1, 0)-\phi_{22}(1)
\right) \alpha(0, t)\nonumber\\
&&+\left( k_{22}(1, 1)+l_{22}(1,1)
\right)\alpha(1, t)
-\sqrt \mu  {\alpha_t}(1, t),\,\,\,\,\label{U_2}
\end{eqnarray}

The main result
is stated next.

\begin{thm}\label{thm1}
Consider system (\ref{plant_eq1})--(\ref{plant_bc2}), with initial conditions $u_0,\alpha_0\in H^1(0,1), u_{0t},\alpha_{0t}\in L^2$,
under the  control law (\ref{U_1})--(\ref{U_2}). If the values of $\delta_1, \delta_2$ (the controller parameters appearing in (\ref{Phi_delta})) are set large enough so that the constant
\begin{equation}
C_2= \min\{\delta_1,\delta_2\}-2-\max\left\{
\frac{4}{\left(\sqrt{\varepsilon}-\theta\right)^2}
,\frac{1}{\mu}\right\},
\end{equation}
is positive, there exists a solution $u(\cdot, t),\alpha(\cdot, t)\in H^1(0,1)$, $u_t(\cdot, t),\alpha_t(\cdot, t)\in L^2(0,1)$ for $t>0$, and the following inequality is verified, guaranteeing the exponential stability of the equilibrium $u\equiv \alpha\equiv u_t \equiv \alpha_t \equiv 0$:
\begin{IEEEeqnarray}{rcl}
&&\|u(\cdot, t)\|^2_{H^1}+\|\alpha(\cdot, t)\|^2_{H^1}+\|u_t(\cdot, t)\|^2_{L^2}+\|\alpha_t(\cdot, t)\|^2_{L^2}
\nonumber\\
&\leq& C_1 \mathrm{e}^{-C_2 t}\Big(\|u_0\|^2_{H^1}+\|\alpha_0\|^2_{H^1}+\|u_{0t}\|^2_{L^2}+\|\alpha_{0t}\|^2_{L^2}
\Big).
\end{IEEEeqnarray}
\end{thm}
The proof of Theorem~\ref{thm1} is given in Section \ref{stability}.

Note that the constant $C_2$ of Theorem~\ref{thm1} only depends on the system parameters $\epsilon$, $\mu$ and $\theta$ and the controller parameters $\delta_1$ and $\delta_2$. Under Assumption~\ref{assump1} it is always possible to set $C_2$ as large as desired, thus achieving arbitrary convergence rate.


\section{Controller Analysis}\label{contdesign}
This section presents the steps leading to (\ref{eqn-controlaw}). The backstepping method is used: first, the target system is presented in Section~\ref{sec-target}; next, the backstepping transformation (of Volterra type) is introduced in Section~\ref{transf}. The well-posedness of the kernel equations is stated in Theorem~\ref{thm2}.
\subsection{Target system}\label{sec-target}
We design a target system as follows
\begin{eqnarray}\label{target10}
\sigma_t&=&\Sigma \sigma_x +\Omega(x) \sigma,
\\ 
\psi_t&=&-\Sigma \psi_x+\Lambda_1 \left( \psi+  \sigma \right)
+\int^x_0 \Xi_2(x, y)\sigma(y, t)dy
\nonumber\\
&&
+\int^x_0 \Xi_3(x, y)\psi(y, t)dy+\Xi_1(x)X\label{target_phi31}, \,\,\,\,\,\\
\dot{X}&=&E_1X+E_2\sigma(0, t),\label{target-X}
\end{eqnarray}
with boundary conditions
\begin{eqnarray}
\sigma(1,t)&=&0,\,
\psi(0,t)=E_3X+C\sigma(0, t),\label{BDtarg}
\end{eqnarray}
where
\begin{eqnarray}\label{eq11}    
\sigma&=&\left[                 
  \begin{array}{c}   
   \eta\\  
    \beta\\ 
  \end{array}
\right],~E_1=(B_1+B_2C)\Phi(0)+A+B_2D,\,\,\,\,\,\,\label{eqn-E1}\\
E_2&=&C\Phi(0)+D,~E_3= B_1+B_2C,
\end{eqnarray}
and where the values of $\Xi_1(x,)$, $\Xi_2(x, y)$, and $\Xi_3(x, y)$ are obtained in terms of the inverse backstepping transformation, in Section~\ref{transf}. The stability of this target system is shown in Section~\ref{stability}.
\subsection{Backstepping transformation}\label{transf}
Firstly, inspired by \cite{c2}, we introduce the following backstepping transformation, of Volterra type
\begin{eqnarray}\label{eq8}
\sigma  &=& Z- \int_0^x {{K}\left( {x,y} \right)Z\left( {y,t} \right)}dy\nonumber
\\\label{tr_eq53}
&&-\int_0^x {{L}\left( {x,y} \right)Y\left( {y,t} \right)}dy-\Phi(x),\\
\psi&=&Y. \label{tr_eq54}
\end{eqnarray}
The kernel equations are deduced as usual, by a tedious but straightforward procedure of taking derivatives in the transformation, replacing the original and target equations, and integrating by parts. The details are skipped for brevity. Regarding their well-posedness, the following result holds. 
\begin{thm}\label{thm2}
 There exists a unique bounded solution
$k_{ij}(x,y), l_{ij}(x, y), i=1,2; j=1,2,3,4$ 
to the kernel equations (\ref{ker_1})--(\ref{ker_5}); in particular, there exists a positive number $M$ such that for $i,j=1,2$
\begin{eqnarray}
|k_{ij}(x,y)|,|l_{ij}(x,y)|&\leq& M\mathrm{e}^{Mx}.
\end{eqnarray}
\end{thm}
The proof follows along the lines of~\cite{c1} and is skipped; it is based on using the method of characteristics to write (\ref{ker_1})--(\ref{ker_5}) in the form of integral equations and then posing a solution in terms of a successive approximation series, whose convergence is proven recursively. It is evident that the derivations of~\cite{c1} can be easily adapted to the presence of integral terms and the  differences in the boundary conditions without much effort.

Since the kernels appearing in (\ref{tr_eq53}) are bounded, the transformation is invertible from the theory of Volterra integral equation. Thus one can define
\begin{eqnarray}
Z &=&\sigma +\int_0^x {{{\mathord{\buildrel{\lower3pt\hbox{$\scriptscriptstyle\smile$}} 
\over K} }}\left( {x,y} \right)\sigma\left( {y,t} \right)} dy\nonumber\\
&&+\int_0^x {{{\mathord{\buildrel{\lower3pt\hbox{$\scriptscriptstyle\smile$}} 
\over L} }}\left( {x,y} \right)\psi\left( {y,t} \right)}dy+{{\mathord{\buildrel{\lower3pt\hbox{$\scriptscriptstyle\smile$}} 
\over \Phi} }}\left( x \right)X,\\\label{eq9}
Y&=&\psi,
\end{eqnarray}
with bounded kernels. Both the transformation and its inverse map $L^2$ functions into $L^2$ functions (see e.g.~\cite{c6}).

From the inverse transformation, the kernels $\Xi_1(x)$, $\Xi_2(x, y)$, $\Xi_3(x, y)$ appearing in (\ref{target_phi31}) are
\begin{eqnarray}
\Xi_1(x)&=&\Lambda_1{{\mathord{\buildrel{\lower3pt\hbox{$\scriptscriptstyle\smile$}} \over \phi} }}(x)+\Lambda_2+\int^x_0 F {{\mathord{\buildrel{\lower3pt\hbox{$\scriptscriptstyle\smile$}} \over \phi} }}(y)dy,\\
\Xi_2(x, y)&=&\Lambda_1 {{\mathord{\buildrel{\lower3pt\hbox{$\scriptscriptstyle\smile$}} \over K} }}(x, y)+F+\int^x_0 F {{\mathord{\buildrel{\lower3pt\hbox{$\scriptscriptstyle\smile$}} \over K} }}(s, y)ds,\\
\Xi_3(x, y)&=&\Lambda_1 {{\mathord{\buildrel{\lower3pt\hbox{$\scriptscriptstyle\smile$}} \over L} }}(x, y)+F+\int^x_0 F {{\mathord{\buildrel{\lower3pt\hbox{$\scriptscriptstyle\smile$}} \over L} }}(s, y)ds,
\end{eqnarray}
from which it can be deduced that they are bounded kernels.

\section{Stability and Well-posedness of Closed Loop}\label{stability}

This section proves Theorem~\ref{thm1}. First, in Section~\ref{sec-solution}, the solution of  (\ref{target10})--(\ref{BDtarg}) is studied with the method of characteristics. This helps to find stability conditions in Section~\ref{stability}. Then, a Lyapunov analysis in Section~\ref{sec-lyap} shows exponential stability. Section~\ref{sec-wp} finishes the proof of Theorem~\ref{thm1} addressing the well-posedness of the system.

\subsection{A semi-explicit solution for the target system}\label{sec-solution}
We start solving (\ref{target10})--(\ref{BDtarg}) with the method of characteristics. Writing down the solution for $\sigma$, it is not difficult to see that 
it converges to zero in finite time $\sqrt{\mu}$. Thus, for $t>\sqrt{\mu}$,
\begin{eqnarray}\label{eq:explicit}
\psi_t(x, t)&=&-\Sigma \psi_x(x, t)+\Lambda_1 \psi(x, t)+\Xi_1(x)X\nonumber\\
&&+\int^x_0 \Xi_3(x, y)\psi(y, t)dy,\nonumber\\
\dot{X}&=&E_1X,\nonumber\\
\psi(0,t)&=&E_3X.
\end{eqnarray}
Solving for $X$ we get $X(t)=X(0)\mathrm{e}^{E_1 t}$, where we have used the matrix exponential. Then
\begin{eqnarray}\label{eq:semi-explicit}
\psi_t(x, t)&=&-\Sigma \psi_x(x, t)+\Lambda_1 \psi(x, t)+\Xi_1(x)X(0)\mathrm{e}^{E_1 t}\nonumber\\
&&+\int^x_0 \Xi_3(x, y)\psi(y, t)dy,\nonumber\\
\psi(0,t)&=&E_3X(0)\mathrm{e}^{E_1 t}.
\end{eqnarray}
Applying the method of characteristics, two Volterra-type integral equations can be found for the components of $\psi$. The details are skipped, but it is easy to see that one can always find an unique $L^2$ solution for $\psi$. 
\subsection{Stability conditions}\label{sub-stability}
Obviously the only requirement for stability is that $E_1$ is Hurwitz as then the origin of the state is exponentially stable for (\ref{eq:explicit}). Nevertheless, for rapid arbitrary stabilization, the eigenvalues of $E_1$ need to be set. Indeed, since $D=-\frac{A}{\sqrt\varepsilon}$ in (\ref{eqn-E1}), $E_1$ is rewritten as
\begin{equation}
E_1=A\left(I-\frac{B_2}{\sqrt\varepsilon}\right) + (B_1+B_2 C) \Phi(0).
\end{equation}
Which, remembering the definition $\kappa=2/(\sqrt{\varepsilon}-\theta)$, results in
\begin{equation}
E_1=\left[
\begin{array}{cc}
 \kappa\xi + \kappa \phi_{11}(0)& -\kappa\left(1+\frac{1}{\sqrt{\mu}}\right)+\kappa \phi_{12}(0)\\
 \frac{\phi_{21}(0)}{\sqrt{\mu}}&\frac{\phi_{22}(0)}{\sqrt{\mu}}\\
\end{array}
\right].
\end{equation}
If we choose the boundary conditions $\Phi(0)$ as follows:
\begin{eqnarray}
\phi_{11}(0)&=&-\xi-\frac{\delta_1}{\kappa}, \,\,
\phi_{12}(0)=1+\frac{1}{\sqrt{\mu}},\label{delta_1_2}\\
\phi_{21}(0)&=&0,\,\,
\phi_{22}(0)=-\delta_2\sqrt{\mu},\label{delta_1_3}
\end{eqnarray}
with $\delta_1,\delta_2>0$, then $E_1$ is a diagonal matrix with entries $-\delta_1$ and $-\delta_2$, which become its eigenvalues. The rate of convergence of $X$ can be arbitrarily set by adjusting the value $\delta_1, \delta_2$ and will be equal to $c=\min \left\{\delta_1,\delta_2\right\}$.

\subsection{Lyapunov-based stability analysis of target system}\label{sec-lyap}
Next, we use a  Lyapunov functional for the stability analysis of target system, to show exponential stability of the origin with a fixed convergence rate. Define
\begin{eqnarray}
V&=&X^T X+\zeta\int^1_0 \mathrm{e}^{\delta x}\sigma^T(x, t) \Sigma^{-1}\sigma(x, t)dx\nonumber\\
&&+ \int^1_0 \mathrm{e}^{-\delta x}\psi^T(x, t) \Sigma^{-1}\psi(x, t)dx.
\end{eqnarray}\label{eq54}
Differentiating (\ref{eq54}) with respect to $t$, we have
\begin{eqnarray}
\dot V&=&2X^T X_t+2\zeta \int^1_0 \mathrm{e}^{\delta x}\sigma^T(x, t)\Sigma^{-1}\sigma_t(x, t)dx\nonumber\\
&&+2\int^1_0 \mathrm{e}^{-\delta x}\psi^T(x, t)\Sigma^{-1}\psi_t(x, t)dx\label{eq59}
\end{eqnarray}
Substituting (\ref{target10})--(\ref{target-X}) into (\ref{eq59}) and then using integration by parts and the fact that $X^T X_t \leq -c X^T X$, we have
\begin{eqnarray}
\dot{V} &\leq&-c X^TX+2X^T E_2\sigma(0, t)-\zeta\sigma^T(0, t)\sigma(0, t)\nonumber\\
&&-\zeta \int^1_0 \mathrm{e}^{\delta x}\sigma^T(x, t) (\delta I+2\Sigma^{-1}\Omega(x))  \sigma(x, t)dx\nonumber
\\
&&- \int^1_0 \mathrm{e}^{-\delta x}\psi^T(x, t)(\delta I-2\Sigma^{-1}\Lambda_1)  \psi(x, t)dx\nonumber\\
&&+2\int^1_0 \mathrm{e}^{-\delta x} \psi^T(x, t)\Sigma^{-1}(\Lambda_1 \sigma(x, t)+\Xi_1 (x)X)dx\nonumber\\
&&+2\int^1_0 \mathrm{e}^{-\delta x} \psi^T(x, t)\Sigma^{-1}\int^x_0 \Xi_2 (x, y)\sigma(y, t)dydx\nonumber\\
&&+2\int^1_0 \mathrm{e}^{-\delta x} \psi^T(x, t)(\Sigma)^{-1}\int^x_0 \Xi_3 (x, y)\psi(y, t)dydx\nonumber\\
&&+ \psi^T(0,t)\psi(0, t)
. \,\,\,\,\,\, \label{eq60}
\end{eqnarray}
Regarding the last line of (\ref{eq60}), using $\psi(0,t)=E_3X+C\sigma(0, t)$, we have
 $ \psi^T(0,t)\psi(0, t)= X^T\mathrm{e}^T_3E_3X+2 X^T\mathrm{e}^T_3C\sigma(0, t)+\sigma^T(0, t)C^TC\sigma(0, t).$
 Then, the first line and last line of (\ref{eq60}) become
\begin{eqnarray}
&&-X^T(c I- \mathrm{e}^T_3E_3)X+2X^T(E_2+ \mathrm{e}^T_3C) \sigma(0, t)\nonumber\\
&&-\sigma^T(0, t)(\zeta- C^TC)\sigma(0, t) \nonumber\\
&\leq & -(c-M_1) X^TX- (\zeta-M_2) \sigma^T(0, t)\sigma(0, t), \label{b66-2}
 \end{eqnarray}
  with $M_1=\max\left\{\kappa^2,\frac{1}{\mu}\right\}+1$, $M_2=\Vert E_2+ \mathrm{e}^T_3C \Vert^2+\Vert C \Vert^2$.
  The fourth line of (\ref{eq60}) is bounded as follows
\begin{eqnarray}
&&2\int^1_0 \mathrm{e}^{-\delta x} \psi^T(x, t)\Sigma^{-1}(\Lambda_1 \sigma(x, t)+\Xi_1 (x)X)dx\nonumber\\
&\leq &
2(1+M_4)\int^1_0 \mathrm{e}^{-\delta x} \psi^T(x, t) \Sigma^{-1} \psi(x,t) dx\nonumber\\
&&+M_3 \int^1_0 \mathrm{e}^{\delta x} \sigma^T(x, t) \Sigma^{-1} \sigma(x,t) dx
+  X^T   X,\label{b67} 
 \end{eqnarray}
 with $M_3=\Vert \Lambda_1 \Vert ^2 $ and $M_4=\max_{x\in[0,1]}\Vert  \Sigma^{-1} \Xi_1 (x) \Vert^2$.
 The fifth line of (\ref{eq60}) can be bounded as follows
\begin{eqnarray}
&&2\int^1_0 \mathrm{e}^{-\delta x} \psi^T(x, t)(\Sigma)^{-1}\int^x_0 \Xi_2(x, y)\sigma(y, t)dydx\nonumber\\
&\leq &
2\int^1_0 \int^1_0 \mathrm{e}^{-\delta x} \vert \psi^T(x, t)\vert \Sigma^{-1}  \vert \Xi_2 (x, y) 
\vert\vert\sigma(y, t)\vert dydx\nonumber\\
&\leq &
 \int^1_0 \mathrm{e}^{-\delta x} \psi^T(x, t) \Sigma^{-1}  \psi(x, t)dx\nonumber\\
&&+M_5 \int^1_0 \mathrm{e}^{\delta x} \sigma^T(x, t)  \Sigma^{-1} \label{b68}
 \sigma(x,t) dx,
 \end{eqnarray}
 where $M_5= \max_{x,y\in[0,1]}\Vert \Xi_2 (x, y) \mathrm{e}^{-\delta y} \Vert^2$.
 Finally, the sixth line of of (\ref{eq60}) is also bounded
\begin{eqnarray}
&&2\int^1_0 \mathrm{e}^{-\delta x} \psi^T(x, t)(\Sigma)^{-1}\int^x_0 \Xi_3 (x, y)\psi(y, t)dydx\nonumber\\
&\leq &
2\int^1_0 \mathrm{e}^{-\delta x/2} \vert \psi^T(x, t)\vert \Sigma^{-1}
\nonumber \\ && \times
\int^x_0  \mathrm{e}^{-\delta x/2}  \vert \Xi_3 (x, y)\vert\vert\psi(y, t)\vert dydx\nonumber\\
&\leq &
2\int^1_0 \int^1_0 \mathrm{e}^{-\delta x/2} \vert \psi^T(x, t)\vert \Sigma^{-1}
\nonumber \\ && \times
\mathrm{e}^{-\delta y/2}  \vert \Xi_3 (x, y) 
\vert\vert\psi(y, t)\vert dydx\nonumber\\
&\leq & M_6
 \int^1_0 \mathrm{e}^{-\delta x} \psi^T(x, t) \Sigma^{-1}  \psi(x, t)dx,
 \label{b69}
 \end{eqnarray}
with $M_6= 1+\max_{x,y\in[0,1]}\Vert \Xi_3 (x, y) \Vert^2$. 
Thus,
\begin{eqnarray}
\dot V&\leq&
 -(c-M_1-1) X^TX- (\zeta-M_2) \sigma^T(0, t)\sigma(0, t) \nonumber
 \\ &&
- \int^1_0 \mathrm{e}^{\delta x}\sigma^T(x, t) (\zeta \delta I+2\zeta \Sigma^{-1}\Omega(x))\sigma(x, t)dx
 \nonumber\\ 
&&
+ \int^1_0 \mathrm{e}^{\delta x}\sigma^T(x, t) (M_3+M_5)\Sigma^{-1}   \sigma(x, t)dx
 \nonumber\\
&&
+ \int^1_0 \mathrm{e}^{-\delta x}\psi^T(x, t)(2\Sigma^{-1}\Lambda_1+M_6\Sigma^{-1}) \psi(x, t)dx\nonumber\\
&&+\int^1_0 \mathrm{e}^{-\delta x}\psi^T(x, t)(3+2M_4-\delta)\psi(x, t)dx.\label{eq70}
\end{eqnarray}
Choosing $c>c'+M_1+1$ with $c'>0$, $\delta >\max \{2 \Vert \Lambda_1 \Vert+M_6+ \Vert \Sigma^{-1} \Vert (c'+1+2M_4),\Vert \Sigma^{-1} \Vert(1+2 \max_{x\in[0,1]} \Vert \Omega(x) \Vert )\}$, and $\zeta >\max\{ M_3+M_5,M_2\}$, we obtain:
\begin{eqnarray}
\dot V&\leq&
 -c' X^TX
-c' \zeta \int^1_0 \mathrm{e}^{\delta x}\sigma^T(x, t)  \Sigma^{-1}\sigma(x, t)dx\nonumber
\\
&&-c' \int^1_0 \mathrm{e}^{-\delta x}\psi^T(x, t)\Sigma^{-1}  \psi(x, t)dx \leq -c' V,
\end{eqnarray}
with $c'>\min \left\{\delta_1,\delta_2\right\}-2-\max\left\{\kappa^2,\frac{1}{\mu}\right\}$. Thus setting the controller parameters $\delta_1$ and $\delta_2$ sufficiently large, an arbitrary convergence rate $c'>0$ is achieved for $V$.

From the Lyapunov inequality just obtained and using norm equivalences, 
and 
the boundedness of the kernels of both direct (\ref{tr_eq53}) and inverse (\ref{eq9}) transformations, 
one obtains 
\begin{IEEEeqnarray}{rcl}
&&\|p(\cdot, t)\|^2_{L^2}+\|q(\cdot, t)\|^2_{L^2}+\|q(\cdot, t)\|^2_{L^2}+\|s(\cdot, t)\|^2_{L^2}
\nonumber\\ &&
+x_1^2(t)+x_2^2(t)
\nonumber\\
&\leq& K_1 \mathrm{e}^{-c' t}\Big(\|p_0\|^2_{L^2}+\|q_0\|^2_{L^2}+\|q_0\|^2_{L^2}+\|s_0\|^2_{L^2}
\nonumber\\ && 
+x_1^2(0)+x_2^2(0)
\Big),
\end{IEEEeqnarray}
for some $K_1>0$. When rewritten in terms of the physical Timoshenko beam states,  the exponential stability bound of Theorem \ref{thm1} follows, since
\begin{eqnarray}
u(t,x)&=&x_1(t)+\frac{1}{2}\int_0^x (p(t,y)+q(t,y)) dy , \\  \label{eqn-uint} 
\alpha(t,x)&=&x_2(t)+\frac{1}{2}\int_0^x (r(t,y)+s(t,y)) dy ,\\   \label{eqn-alphaint} 
 u_t(t,x)&=&\frac{p(t,x)-q(t,x)}{2\sqrt\epsilon},\\
 \alpha_t(t,x)&=&\frac{r(t,x)-s(t,x)}{2\sqrt\mu}.\label{eqn-alphatint} 
 \end{eqnarray}
\subsection{Well-posedness of the closed-loop system}\label{sec-wp}
Under the assumptions of Theorem 1, we have that $u_0\in H^1, \alpha_0\in H^1, u_{0t}\in L^2, \alpha_{0t}\in L^2$
and therefore the initial conditions of $p,q,r,s$ belong to $L^2$. Therefore the initial conditions of the transformed states are also $L^2$. It is easy to see that the target system is well-posed in $L^2$ (see Section IV.B.1); thus the original system will be as well, since the inverse  transformation maps $L^2$ into $L^2$. This finally implies the well-posedness result of Theorem 1, by  (\ref{eqn-uint})--(\ref{eqn-alphatint}).
\section{Numerical Simulation}\label{sec-numerical}
\begin{figure*}[t!]
      \begin{center}
      \includegraphics[width=5.8cm]{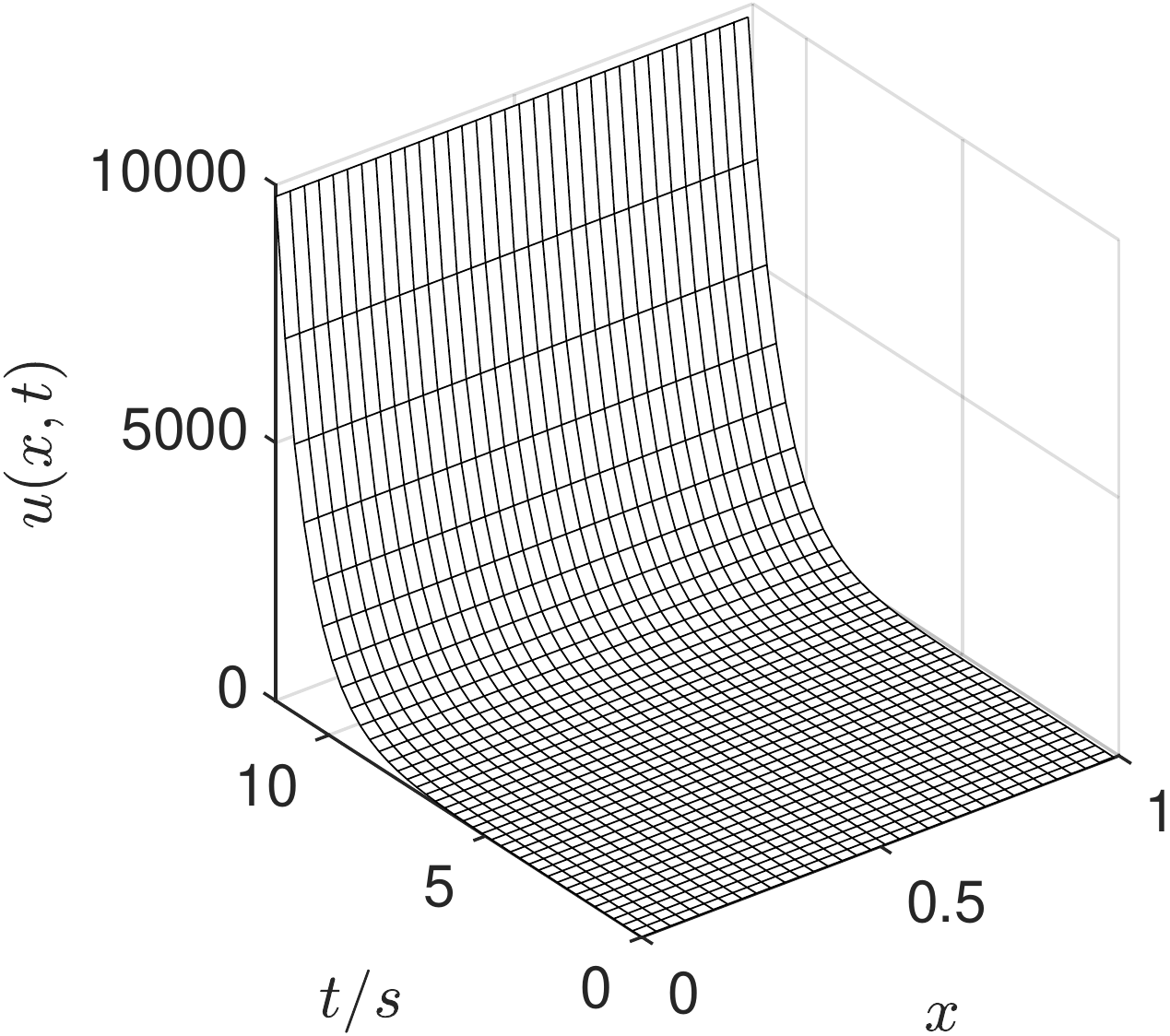}
      \includegraphics[width=5.8cm]{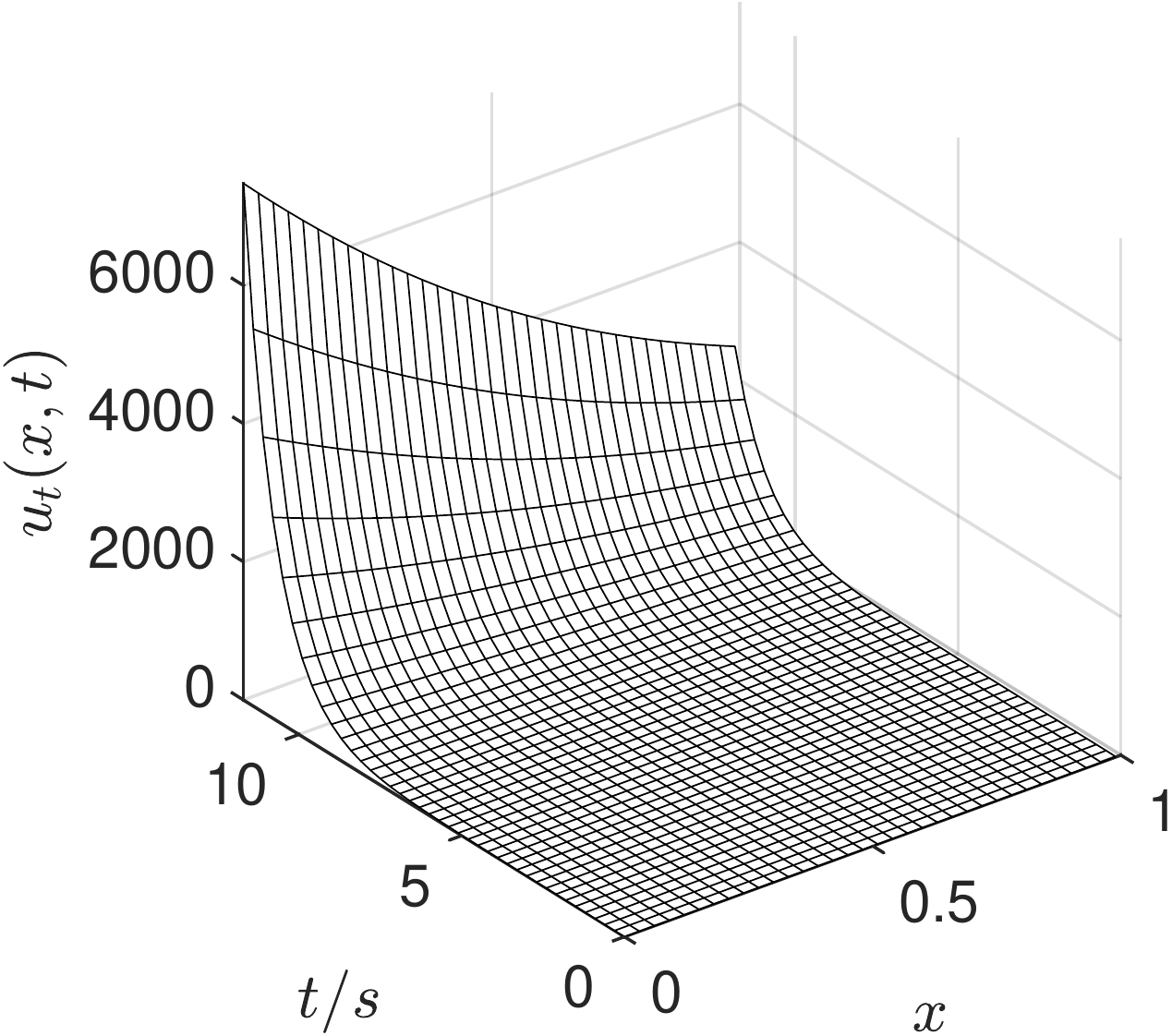}\\
      \includegraphics[width=5.8cm]{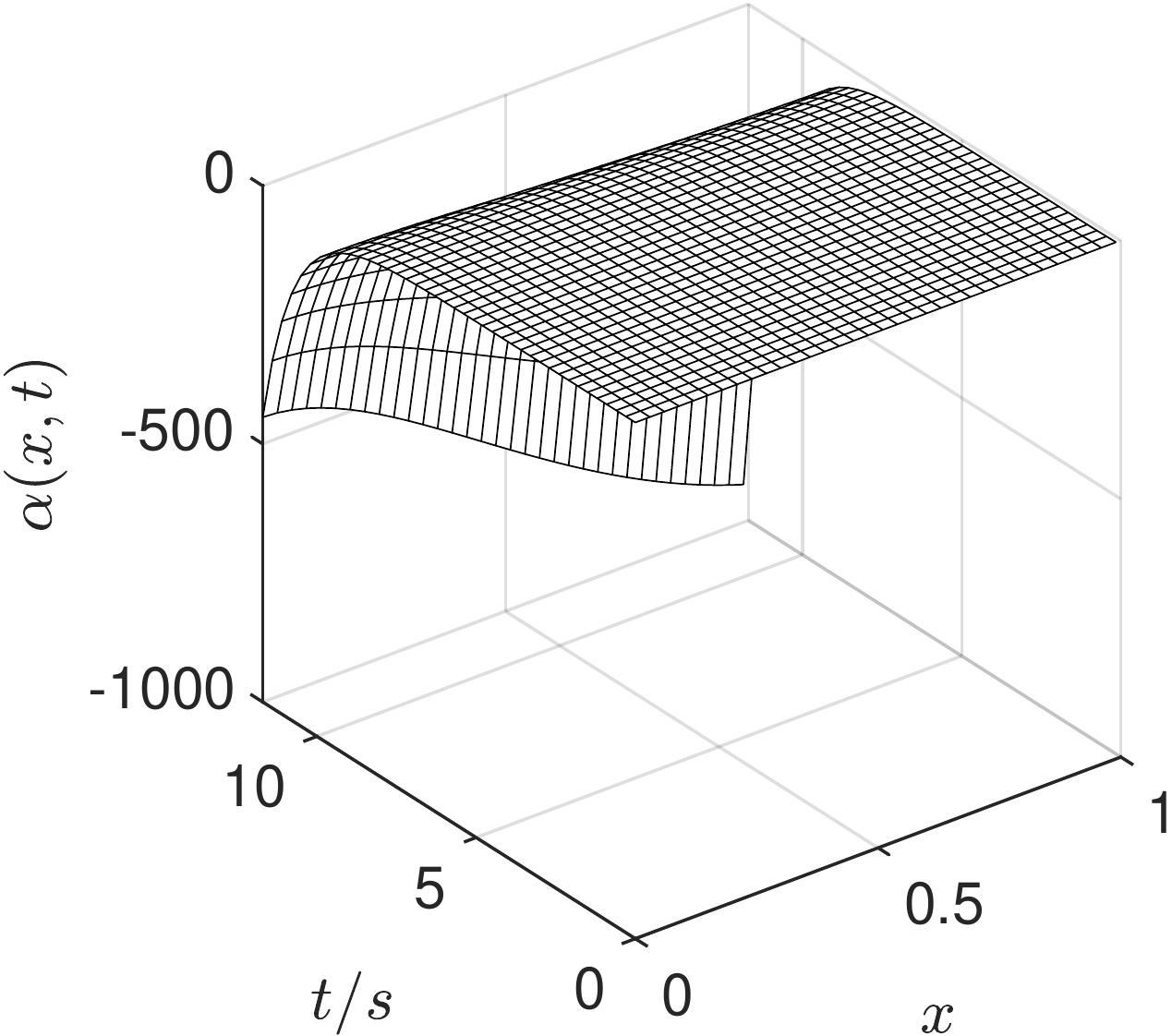}
      \includegraphics[width=5.8cm]{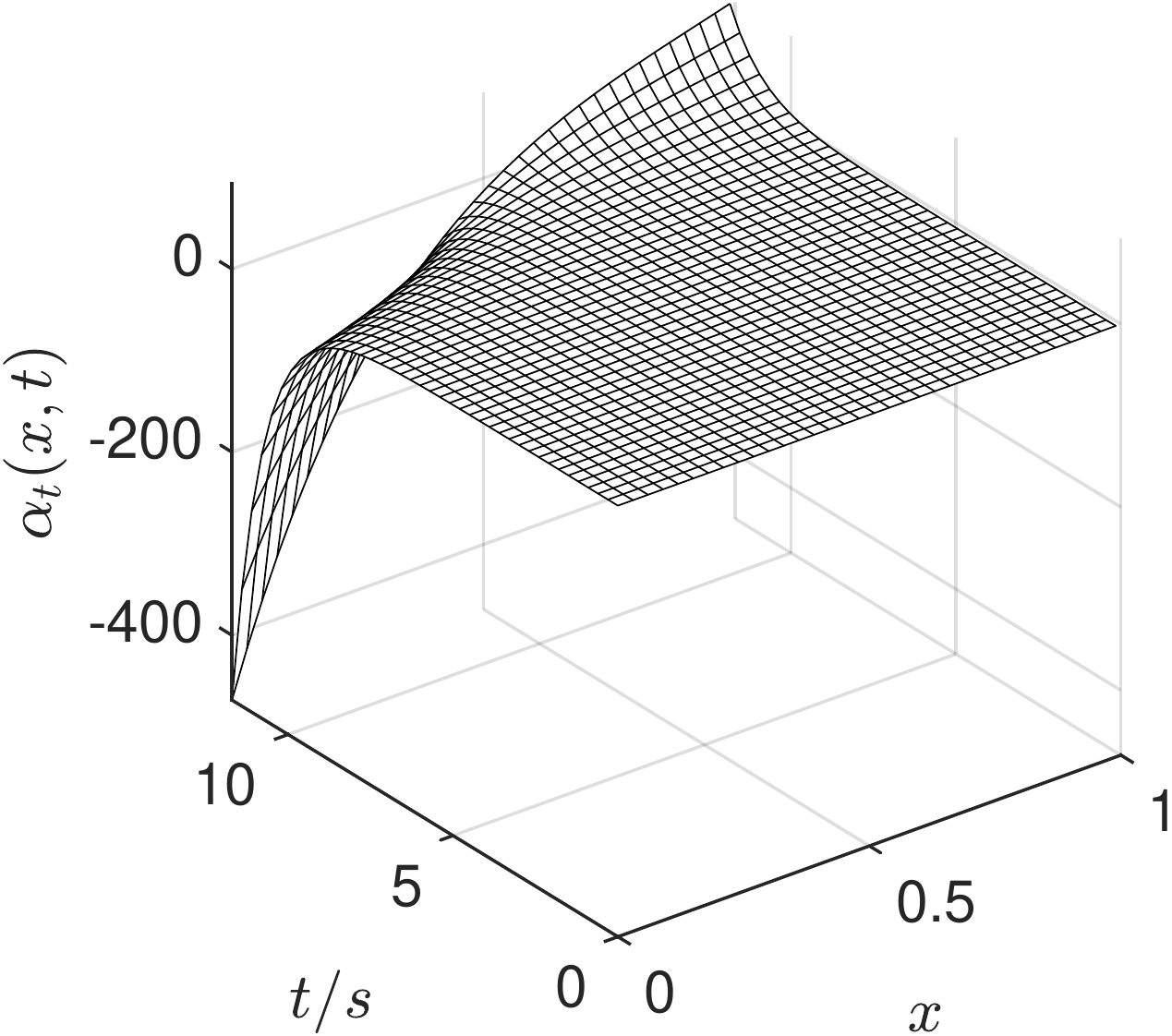}
      \caption{Open-loop evolution of Timoshenko beam states $u$, $u_t$, $\alpha$, $\alpha_t$ over time.}
      \label{Fig-openloop}
      \end{center}
\end{figure*}
\begin{figure*}[t!]
     \begin{center}
      \includegraphics[width=5.7cm]{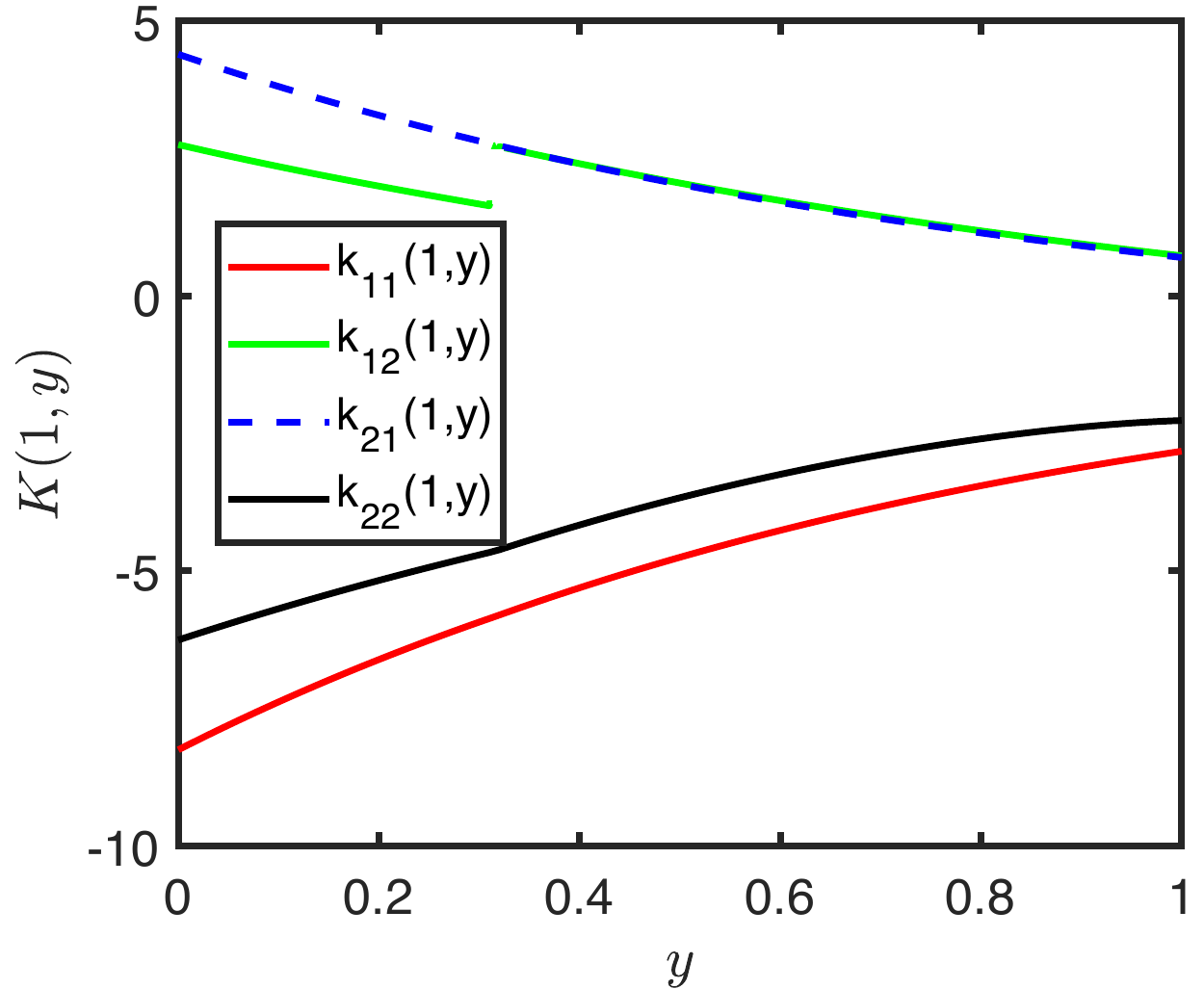}
      \includegraphics[width=5.7cm]{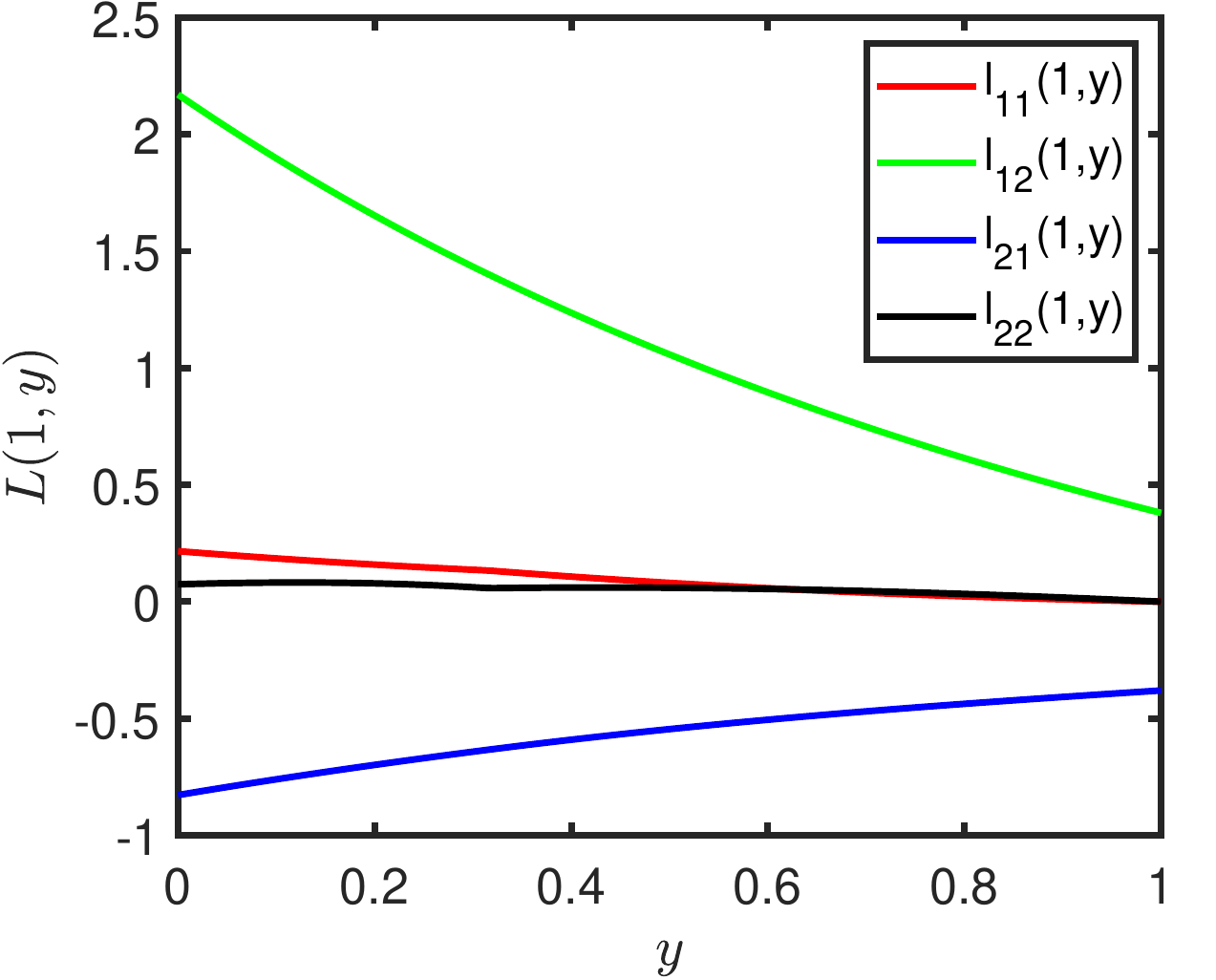}
      \includegraphics[width=5.7cm]{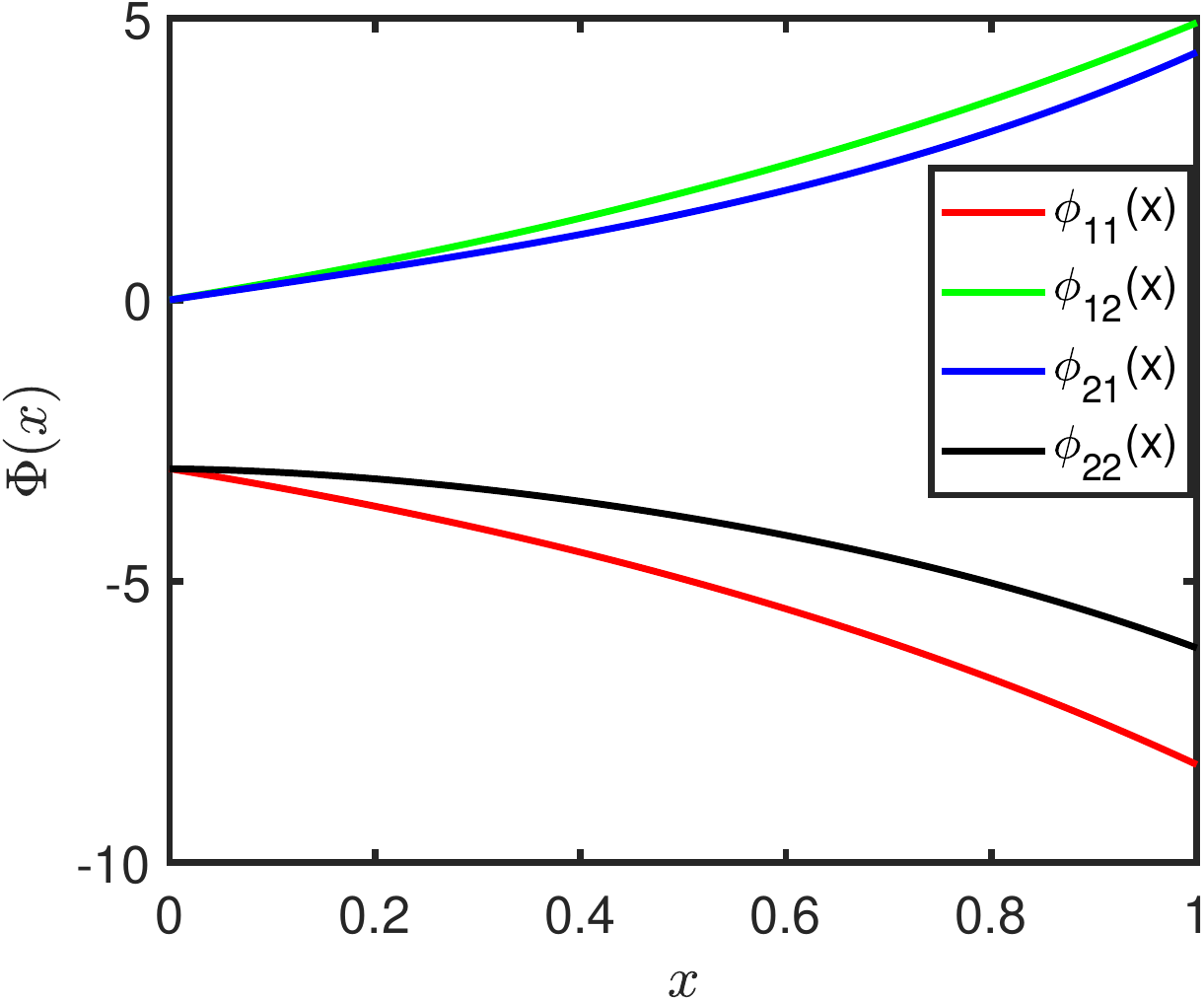}
      \caption{Feedback gains $K(1,y)$, $L(1,y)$, $\Phi(x)$. Note the discontinuity in the kernel function $k_{12}(1,y)$.}
      \label{Fig-gains}
      \end{center}
 \end{figure*}

\begin{figure*}[t!]
      \begin{center}
      \includegraphics[width=5.8cm]{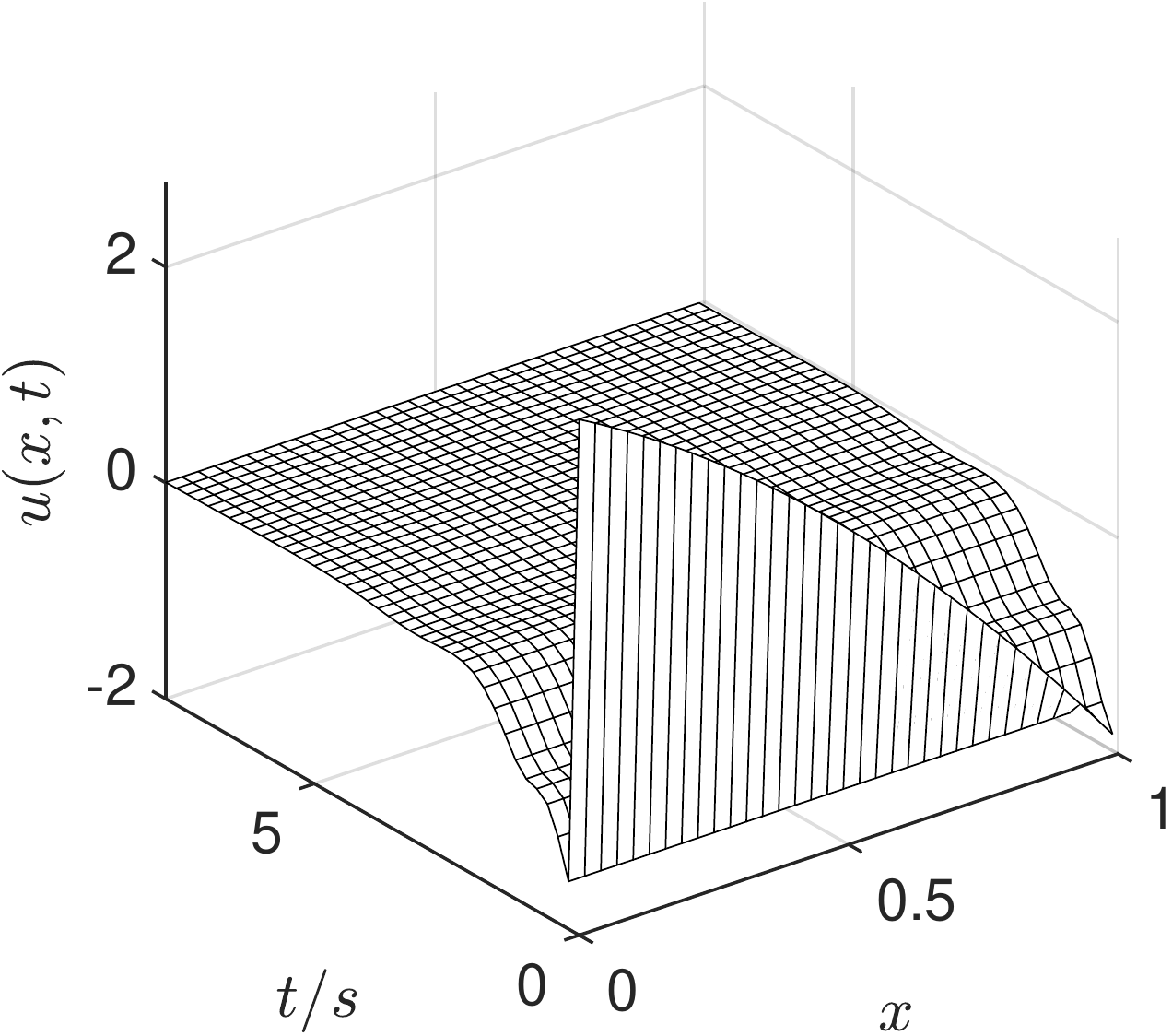}
      \includegraphics[width=5.8cm]{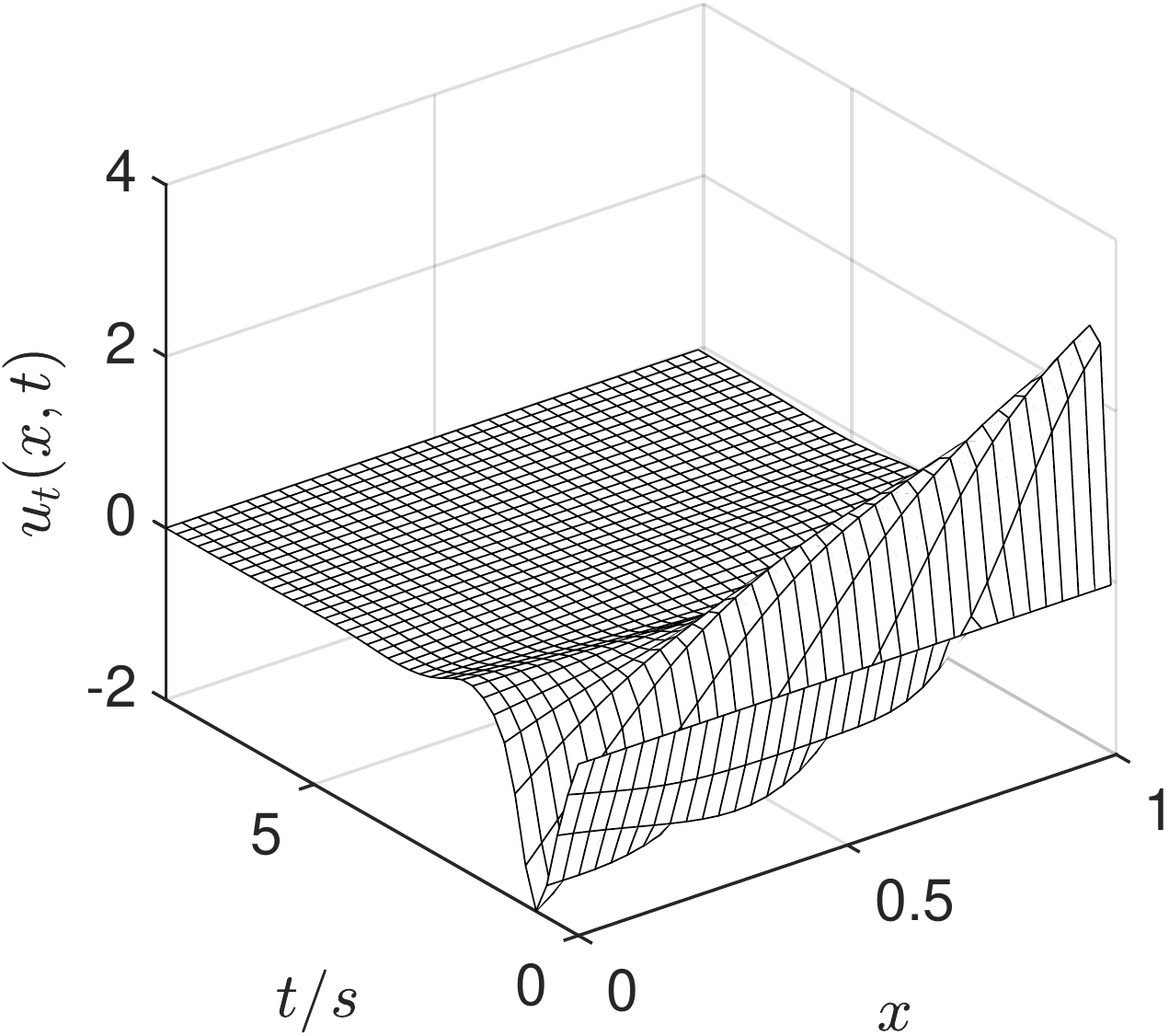}\\
      \includegraphics[width=5.8cm]{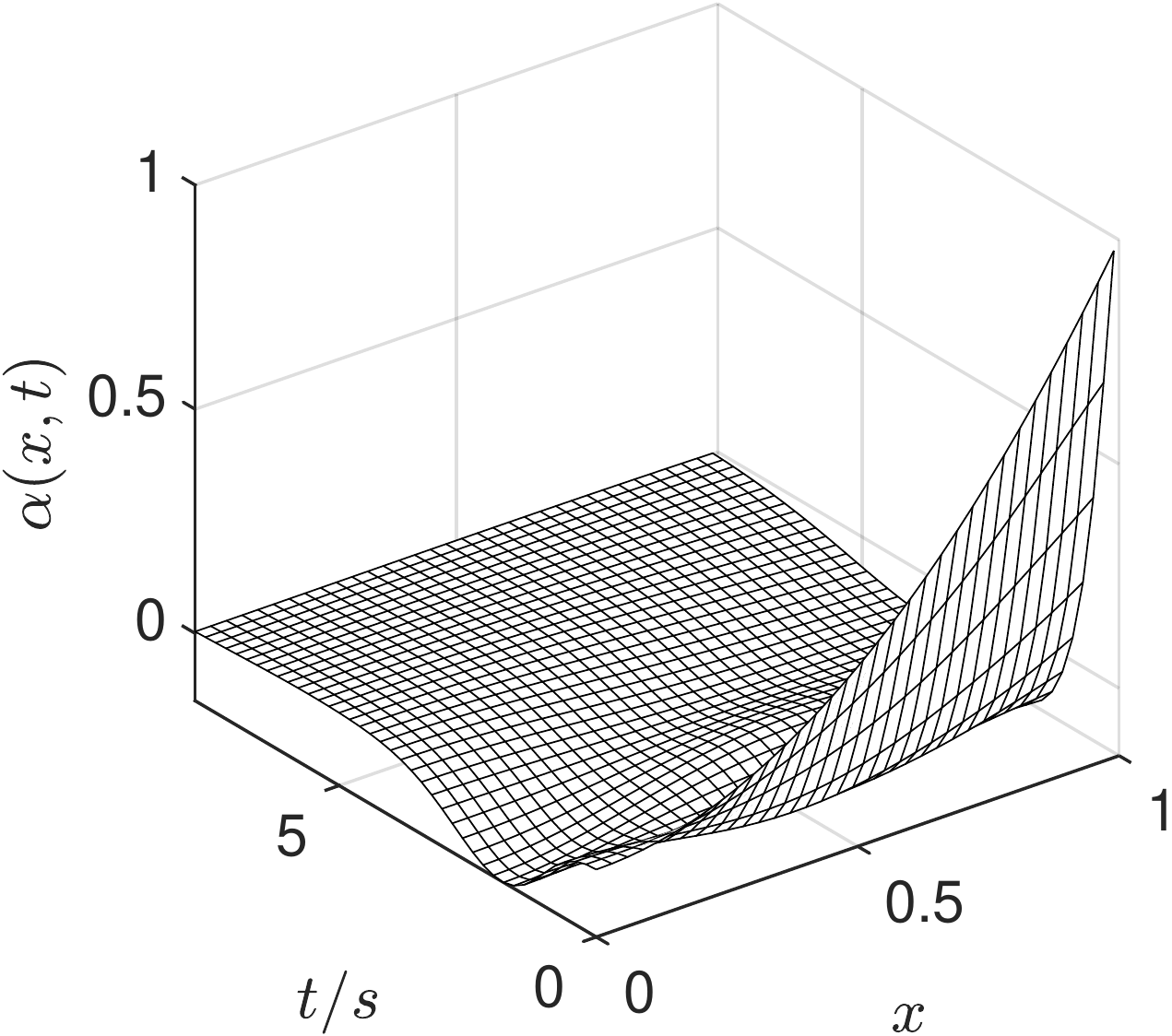}
      \includegraphics[width=5.8cm]{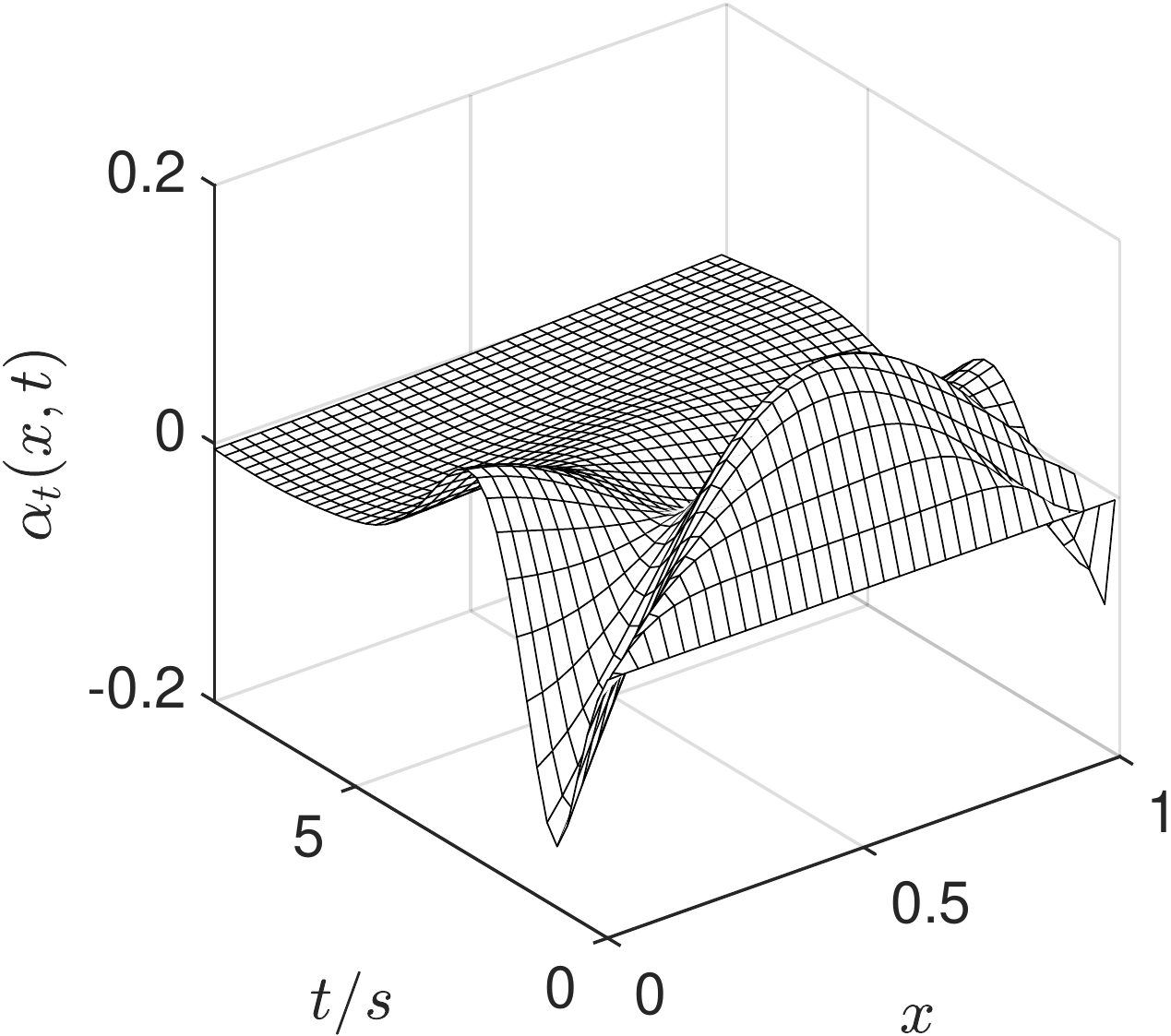}
      \caption{Closed-loop evolution of Timoshenko beam states $u$, $u_t$, $\alpha$, $\alpha_t$ over time.}
      \label{Fig-closedloop}
      \end{center}
          \end{figure*}

 
To verify the effectiveness of the proposed boundary controller, (\ref{plant_eq1})--(\ref{plant_bc2}) is simulated with $\epsilon=1$, $\mu=2$, $a=1$, $\theta=-1$, $\xi=1$. The initial values are set to  $u_{0}=2.8-2.8x-1.8x^2, u_{t0}=0, \alpha_{0}=x^2, \alpha_{t0}=0$.
We use the HPDE tool in MATLAB, in which the four equivalent one-order hyperbolic PDEs (\ref{Equ_eq1})--(\ref{Equ_eq4}) and the ODEs (\ref{ODE_eq1})--(\ref{ODE_eq2}) are solved,  and the evolution of $u(x, t), \alpha(x, t)$ is obtained by using (\ref{eqn-uint})--(\ref{eqn-alphaint}). We first show in Fig.~\ref{Fig-openloop} the unstable response of the open-loop system, which diverges due to anti-damping. Next, we apply the proposed controller (\ref{U_1})--(\ref{U_2}) to the Timoshenko beam. The controller parameters are chosen as $\delta_1=5, \delta_2=2$. The feedback gains $K(1,y)$, $L(x,y)$ and $\Phi(x)$ are shown in Fig.~\ref{Fig-gains} and were computed using a power series approach as in \cite{Leo2020}. There is a discontinuity in the kernel function $k_{12}(1,y)$ , which is typically present when applying the backstepping method to a $(2+2)\times(2+2)$ system and does not impact the result~\cite{c6}. The variables $u(x, t)$, $u_t(x, t)$, $\alpha(x, t)$ and $\alpha_t(x, t)$ evolve as shown in Fig.~\ref{Fig-closedloop}, converging to zero exponentially, as expected from Theorem~\ref{thm1}. 
\section{Concluding remarks}\label{sec-conclusion}
This work considered boundary control of a Timoshenko beam with anti-damping and anti-stiffness at the uncontrolled boundary; firstly, we transform the Timoshenko beam states into a hyperbolic PIDE-ODE system. Then, backstepping is applied, obtaining arbitrarily fast decay. Simulations show the effectiveness of the  controller in setting the convergence rate. As future work, the possibility of achieving \emph{finite-time} convergence (in the spirit of the superstability results of \cite{bal1997}), by using time-varying backstepping~\cite{steeves2021}, will be investigated.


\begin{thebibliography}{99}
\bibitem{Ammar2007}
Ammar-Khodja, F., Kerbal, S., and Soufyane, A., ``Stabilization of the nonuniform Timoshenko beam,'' {\it Journal of Mathematical Analysis and Applications,} 327(1), 525--538, 2007.
\bibitem{bal1997}
Balakrishnan, A.V., ``On Superstability of semigroups,'' In {\it Proceedings of the 18th IFIP TC-7}, pp. 12--19, 1997.
\bibitem{bal2005}
Balakrishnan, A. V., ``Superstability of systems,'' {\it Appl. Math. Comput.}, 164, pp. 321--326, 2005.
\bibitem{c5}
Bastin, G. , and Coron, J.-M., {\it Stability and Boundary Stabilization of 1-D Hyperbolic Systems}, Basel: Birkhäuser, 2016.
\bibitem{Leo2020}
Camacho-Solorio, L., Vazquez R., and Krstic, M.,``Boundary observers for coupled diffusion–reaction systems with prescribed convergence rate,'' {\it Systems \& Control Letters} 135 : 104586, 2020.
\bibitem{chen2022}
Chen, G., Vazquez, R., and Krstic, M.,`` Backstepping-Based Exponential Stabilization of Timoshenko Beam with Prescribed Decay Rate,'' submitted to {\it IFAC CPDE 2022}, 2022.
\bibitem{Coron1998}
Coron, J.M. and d'Andrea-Novel, B., ``Stabilization of a rotating body beam without damping,'' {\it IEEE T. Automat. Contr.}, 43(5), pp. 608--618, 1998.
\bibitem{c1}
Di Meglio, F., Argomedo, F. B. , Hu, L., and  Krstic, M., ``Stabilization of coupled linear heterodirectional hyperbolic PDE-ODE systems,'' {\it Automatica},
 vol. 87, pp. 281-289, 2018.
 \bibitem{Guo2021}
 Guo, B.Z. and Meng, T., ``Robust output regulation for Timoshenko beam equation with two inputs and two outputs,'' {\it International Journal of Robust and Nonlinear Control}, 31(4), pp.1245--1269, 2021. 
\bibitem{Han1999}
  Han, S. M., Benaroya, H.,  and Wei, T., “Dynamics of transversely vibrating beams using four engineering theories,” {\it Journal of Sound and Vibration}, vol. 225, pp. 935-–988, 1999.
 \bibitem{Wei2015}
He, W., and Liu, C., ``Vibration control of a Timoshenko beam system with input backlash,'' {\it  IET Control. Theory Appl.,} 9(12), 1802--1809, 2015.
 \bibitem{HeZhang2012}
He, W., Zhang, S., and Ge, S. S., ``Boundary output-feedback stabilization of a Timoshenko beam using disturbance observer,'' {\it IEEE Transactions on Industrial Electronics,} 60(11), 5186--5194, 2012.
 \bibitem{He2014}
He, W., Zhang, S., Ge, S. S., and Liu, C., ``Adaptive boundary control for a class of inhomogeneous Timoshenko beam equations with constraints,'' {\it IET Control Theory \& Applications}, 8(14), 1285--1292, 2014.
\bibitem{c6}
Hu, L., Vazquez, R., Di Meglio, F., Krstic, M., "Boundary exponential stabilization of 1-D inhomogeneous quasilinear hyperbolic systems," {\it SIAM J. Contr. Optim.}, vol. 57, pp. 963--998, 2019.
 \bibitem{Kim1987}
Kim, J. U., and Renardy, Y., ``Boundary control of the Timoshenko beam,'' {\it  SIAM J. Contr. Optim.}, 25(6), 1417--1429, 1987.
\bibitem{Krstic2006}
Krstic, M., Siranosian, A. A., and Smyshlyaev, A., ``Backstepping boundary controllers and observers for the slender Timoshenko beam: Part I-Design,'' {\it In 2006 ACC}, pp. 2412--2417, 2006.
\bibitem{Krstic2008} Krstic, M., Guo, B.Z., Balogh, A. and Smyshlyaev, A., ``Control of a tip-force destabilized shear beam by observer-based boundary feedback,'' {\it  SIAM J. Contr. Optim.}, 47(2), pp.553--574, 2008.
\bibitem{c2}
Lingling, S., Wang, J.-M., and Krstic, M., ``Boundary feedback stabilization of a class of coupled hyperbolic equations with nonlocal terms,'' {\it IEEE T. Automat. Contr.}, vol. 63, no. 8, pp. 2633-2640, 2017.
\bibitem{Liu2013}
Liu, X. F., and Xu, G. Q., ``Exponential stabilization for Timoshenko beam with distributed delay in the boundary control,'' {\it In Abstract and Applied Analysis},  Hindawi, 2013.
\bibitem{Ma2020}
Ma, J., Wei, Z., Wen, H., and Jin, D., ``Boundary control of a Timoshenko beam with prescribed performance,'' {\it Acta Mechanica,} 231, 3219--3234, 2020.
\bibitem{Macchelli2004}
Macchelli, A., and Melchiorri, C., ``Modeling and control of the Timoshenko beam: The distributed port Hamiltonian approach,'' {\it SIAM J. Contr. Optim.,} 43(2), 743--767, 2004.
\bibitem{mattioni2020}
Mattioni, A., Wu, Y., Le Gorrec, Y. and Zwart, H., ``Stabilisation of a rotating beam clamped on a moving inertia with strong dissipation feedback,'' in {\it 59th IEEE CDC}, pp. 5056--5061, 2020.
\bibitem{mattioni2021}
Mattioni, A., Wu, Y. and Le Gorrec, Y., ``Exponential stabilization of a clamped Timoshenko beam with actuation on a tip mass,'' in {\it 60th IEEE CDC}, pp. 6200--6205, 2021.
\bibitem{Mor1992}
Morgül, O., ``Dynamic boundary control of the Timoshenko beam,'' {\it Automatica,} 28(6), 1255--1260, 1992.
\bibitem{Siuka2011}
Siuka, A., Schöberl, M., and Schlacher, K., ``Port-Hamiltonian modelling and energy-based control of the Timoshenko beam,'' {\it Acta mechanica,} 222(1), 69--89, 2011.
\bibitem{Siranosian2009}
Siranosian, A. A., Krstic, M., Smyshlyaev, A., and Bement, M., ``Motion planning and tracking for tip displacement and deflection angle for flexible beams,'' {\it J. Dyn. Sys. Meas. Control,} 131(3), 2009
\bibitem{Smysh2009} Smyshlyaev, A., Guo, B.Z. and Krstic, M., ``Arbitrary decay rate for Euler-Bernoulli beam by backstepping boundary feedback,'' {\it IEEE Transactions on Automatic Control}, 54(5), pp.1134--1140, 2009.
\bibitem{steeves2021}
Steeves, D., and Krstic, M., ``Prescribed-time stabilization of ODEs with diffusive actuator dynamics,'' in {\it 24th MTNS}, 2021.
\bibitem{Soufyane2003}
Soufyane, A., and Wehbe, A., ``Uniform stabilization for the Timoshenko beam by a locally distributed damping,'' {\it Electron. J. Differ. Equ.}, 29, pp. 1--14, 2003.
\bibitem{Taylor2002}
Taylor, S. W., and Yau, S. C., `` Boundary control of a rotating Timoshenko beam,'' {\it ANZIAM Journal}, 44, E143--E184,2002.
\bibitem{Wang2015} Wang, J.M. and Krstic, M., ``Stability of an interconnected system of Euler-Bernouilli beam and heat equation with boundary coupling,'' {\it ESAIM: COCOV}, 21(4), pp.1029--1052, 2015.
\bibitem{Wu2017}
Wu, Y., Hamroun, B., Le Gorrec, Y. and Maschke, B., ``Reduced order controller design for Timoshenko beam: A port Hamiltonian approach,'' {\it IFAC-PapersOnLine}, 50(1), pp.7121--7126, 2017.
\bibitem{Xu2003}
Xu, G. Q., and Yung, S. P., ``Stabilization of Timoshenko beam by means of pointwise controls,'' {\it ESAIM: COCOV}, 9, 579--600, 2003.
\bibitem{Xu2005}
Xu G. Q., ``Boundary feedback exponential stabilization of a Timoshenko beam with both ends free,'' {\it International Journal of Control,} 78(4), 286--297, 2005.
\bibitem{Yildirim2016}
Yildirim, K., and Kucuk, I., ``Active piezoelectric vibration control for a Timoshenko beam,'' {\it  J. Frankl. Inst.,} 353(1), 95--107, 2016.
\bibitem{Zhao2005}
Zhao, H. L., Liu, K. S., and Zhang, C. G., ``Stability for the Timoshenko beam system with local Kelvin–Voigt damping,'' {\it Acta Mathematica Sinica,} 21(3), 655--666, 2005.











\end{thebibliography}
\end{document}